\documentclass{amsart}

\usepackage{amsfonts,amsmath,upref,amssymb,eucal}

\newtheorem{theorem}{Theorem}
\newtheorem{lemma}[theorem]{Lemma}

\newenvironment{definition}{\medskip \refstepcounter{theorem}
\noindent  {\bf Definition \thetheorem}.\rm}{\,}

\newenvironment{remark}{\medskip \refstepcounter{theorem}
\noindent  {\bf Remark \thetheorem}.\rm}{\,}

\def\a{\alpha}
\def\b{\beta}
\def\d{\delta}
\def\qp{{\mathbb Q}_p}
\def\zp{{\mathbb Z}_p}

\def\z{\mathbb Z}

\def\np#1{\left| \! \left| #1 \right| \! \right|_p }

\def\ls#1#2{\left( \frac{#1}{#2} \right)}

\begin{document}
\title{Arithmetic differential operators on $\zp$}
\author[A. Buium]{Alexandru Buium}
\thanks{During the preparation of this work, the first author
was partially supported by NSF grant DMS 0552314.}
\author[C.C. Ralph]{Claire C. Ralph}
\author[S.R. Simanca]{Santiago R. Simanca}
\address{Department of Mathematics \& Statistics, The University of New
Mexico, NM 87131}
\email{buium@math.unm.edu, cralph@unm.edu, santiago@math.unm.edu}

\begin{abstract}
Given a prime $p$, we let $\d x=(x-x^p)/p$ be the the Fermat quotient
operator over $\zp$. We prove that a function
$f:\zp \rightarrow \zp$ is analytic if, and only if,
there exists $m$ such that $f$ can be
represented as $f(x)=F(x, \d x, \ldots , \d^m x)$, where $F$ is a
restricted power series with $\zp$-coefficients in $m+1$ variables.
\end{abstract}

\maketitle

\section{Introduction}
\subsection{Main results}
An arithmetic analogue of the theory of ordinary differential equations
was initiated in \cite{ab3}, and further developed in a series of subsequent
publications; see \cite{ab2}, and the bibliography therein.
 In this theory, the role of the differentiation operator is played by a
Fermat quotient operator acting on numbers.
Later on, the theory was extended to one of arithmetic partial
differential equations \cite{alsa,alsa1,alsa2}. Our work here derives its
motivation from certain examples encountered in the original development
of the theory, but is independent of them.

Let $p$ be a prime integer that we fix hereafter. We denote by
$\qp$ the field of $p$-adic numbers, with $p$-adic norm $\np{\, \cdot \,}$
normalized by $\np{p}=p^{-1}$. On the ring
$\zp=\{x\in \qp: \np{x}\leq 1\}$
of $p$-adic integers, we consider the {\it Fermat quotient operator}
\begin{equation}
\d a= \frac{a-a^p}{p}\, ,\label{do}
\end{equation}
which by analogy we view as the ``derivative of $a$ with respect to $p$.''
We denote by $\d^i $ the $i$-th iterate of $\d$.

Given a multi-index $\a=(\a_0, \ldots, \a_k)$ of non-negative integers, we
shall say that $\a \geq 0$, and use the expression $x^{\a}$ to
denote the monomial $x_0^{\a_0}\cdots x_k^{\a_k}$. By $|\a|$
we mean $|\a|= \a_0 + \cdots + \a_k$. We recall that
$F(x)=\sum_{\a\geq 0}a_{\a}x^{\a} \in \zp[[x_0, \ldots, x_k]]$ is said to
be a {\it restricted} power series if $\lim_{|\a| \rightarrow \infty}a_{\a}=0$.

\begin{definition}
\label{firstdef}
A function $f:\zp \rightarrow \zp$ is called an
{\it arithmetic differential operator of order} $m$
if there exists a restricted power series $F \in \zp[[x_0,x_1,\ldots ,x_m]]$
such that
\begin{equation}
\label{bubu}
f(a)=F(a,\d a, \ldots ,\d^m a)
\end{equation}
for all $a \in \zp$. We say that the series $F$ {\it $\d$-represents} $f$.
\qed
\end{definition}

\begin{definition}
\label{def2}
A function $f:\zp \rightarrow \zp$ is said to be {\it analytic of level} $m$,
if for any $a \in \zp$ there exists a restricted power series
$F_a\in \zp[[x]]$ such that
$$
f(a+p^m u)=F_a(u)
$$
for  all $u \in \zp$. We say that the collection of series
$F_a$ {\it represents} $f$.
\qed
\end{definition}

\begin{remark}
\label{trei}
We have the following simple observations:
\begin{enumerate}\label{r1}
\item In Definition \ref{def2}, it is enough to take the $a$s in a
complete residue system mod $p^m$ in $\zp$.
\item A function $f: \zp \rightarrow \zp$ is analytic in the sense of
\cite{ma} (cf. with \cite{se}, p. LG 2.4) if, and only if,
it is analytic of level $m$ for some $m$.
\item If $f,g:\zp \rightarrow \zp$ are arithmetic differential operators
of orders $m$ and $n$, respectively, then $f \circ g$ is an arithmetic
differential operator of order $m+n$.
\item Any arithmetic differential operator $\zp \rightarrow \zp$ is an
analytic function in the sense of \cite{ma,se}. \qed
\end{enumerate}
\end{remark}
\medskip

It will be relatively easy to sharpen the last part of this remark to the
following:

\begin{theorem}
\label{easy}
Any arithmetic differential operator  $f: \zp \rightarrow \zp$ of order $m$ is
an analytic function of level $m$.
\end{theorem}

Our main result says that the converse of this statement is true:

\begin{theorem}
\label{main}
Any analytic function $f: \zp \rightarrow \zp$ of level $m$ is an
arithmetic differential operator of order $m$.
\end{theorem}

\begin{remark}
Actually, we will also prove that among the restricted power series
that $\d$-represent $f$, there is a unique one
$F=F(x_0,\ldots ,x_m)$ satisfying the condition that all of its monomials
are of degree less or equal that $p-1$ in each of the variables
$x_0,\ldots ,x_{m-1}$; cf.
Theorem \ref{mainn}.
We call this $F$ the {\it canonical} series $\d$-representing $f$.
As the proof of Theorem \ref{mainn} will show, the computation of the
canonical series $F$ $\d$-representing $f$ in terms of the collection of
series $F_a$ representing $f$ is a rather non-trivial task.
\end{remark}

\begin{remark}
The context of our work here and that of the theory in \cite{ab2,ab3}
differ from each other. Indeed, in \cite{ab2,ab3} the ring $\zp$ is replaced
by the completion $R$ of the maximum unramified extension of $\zp$,
and our $\d$ here is replaced by
$$
\begin{array}{ccc}
R & \stackrel{\d}{\rightarrow} & R \\
a & \mapsto & {\displaystyle \frac{\phi(a)-a^p}{p}}\, ,
\end{array}
$$
where $\phi$ is the unique lift of Frobenius on $R$. In this other
setting, an arithmetic differential operator that is not of order $0$ is
never analytic, and an analytic function that is not of level $0$
is never an arithmetic differential operator. The difference between the
theories over these two rings is, in some sense, analogous to the difference
between number theoretic statements about finite fields and algebraic
geometric statements over their algebraic closures.
\qed
\end{remark}

\begin{remark}
If we use a slightly more general context, the
theory in \cite{ab2, ab3} gives rise to several interesting
number theoretic locally constant functions that have nice
representations as arithmetic differential operators of low order.
As these have been the main source of motivation for our work, we
describe briefly two of them here.

Definition \ref{firstdef} can be extended by considering
functions $f: \zp^N \rightarrow \zp$ that can be represented as
in (\ref{bubu}) with
an $F\in \zp[[x_0, \ldots, x_m]]$ but where now each $x_j$ is an
$N$-tuple of variables. We call these {\it arithmetic differential operators
of order $m$} also. If $X$ is an affine scheme embedded into the affine
$N$-space over $\zp$,
we let $X(\zp) \subset \zp^N$ be the natural inclusion at the level of
$\zp$-points. Then a function $X(\zp) \rightarrow \zp$ is called an
{\it arithmetic differential
 of order} $m$ if it can be extended to an arithmetic differential operator
$\zp^N \rightarrow \zp$ of order $m$.

For instance, let $X$ be the multiplicative group scheme over $\zp$
embedded into the affine plane ${\rm Spec}\, \zp [v,w]$ via the map
$u \mapsto (u,u^{-1})$. The we have that
$X(\zp)=\zp^{\times}$, and we can
talk about arithmetic differential operators
$\zp^{\times} \rightarrow \zp$  of order
$m$. As noted in \cite{ab1}, for odd primes $p$,
the {\it Legendre symbol} $f:\zp^{\times} \rightarrow \zp$, defined by
$$f(a):=\ls{a}{p}=\left\{ \begin{array}{rl}
1  &  \text{if $a$ (mod $p$) is a quadratic residue mod $p$}, \\
-1 &  \text{if $a$ (mod $p$) is a quadratic nonresidue mod $p$},
\end{array} \right.
$$
is the arithmetic differential operator of order one given by
$$
\ls{a}{p} = a^{\frac{p-1}{2}}\left( 1+ \sum_{n=1}^{\infty}
(-1)^{n-1}\frac{(2n-2)! p^n }{2^{2n-1}(n-1)!n!} (\d a)^n
(a^{-1})^{-pn}\right)\, .
$$
This function is locally constant of level $1$, that is to say, constant on
discs of radius $1/p$.

Similarly, let $X$ be the locus in the plane
${\rm Spec}\, \zp[u,v]$ over $\zp$, where $4v^3+27w^2$ is invertible,
and view $X$ as embedded in $3$-space via the map $(v,w)\mapsto
(v,w,(4v^3+27w^2)^{-1})$. Consider the traces of
Frobenii of the reductions mod $p$ of elliptic curves $y^2=x^3+Ax+B$ over
$\zp$. These functions can be represented as quotients of some remarkable
arithmetic
differential operators of order $2$ defined on
$$
X(\zp)=\{(A,B)\in \zp\times\zp: \,
4A^3+27B^2\in \zp^{\times}\}\, .
$$
Cf. \cite{ab1} for details.
\qed
\end{remark}

\subsection{Continuous functions on $\zp$}
It is of interest to compare our results with a known Mahler-type theorem
about the structure of continuous $\zp$-valued functions on
$\zp$.

Let $A$ be the set of all non-negative integral
vectors $\a=(\a_0, \a_1, \a_2,\ldots )$ with finite support.
Thus, $\a_j \geq 0$ for all $j$, and
$\a_j=0$ for $j$ sufficiently large. If $\a \in A$, it makes sense
to compute $|\a|$. Given a sequence of variables $x_0,x_1,x_2,\ldots$,
we set $x^{\a}$ for $x_0^{\a_0} x_1^{\a_1} x_2^{\a_2}\ldots $.
Then we say that a power series
$$
F(x_0,x_1,x_2,\ldots )=\sum_{\a\in A} a_{\a} x^{\a}\, , \; a_{\a}\in \zp\, ,
$$
is {\it restricted} if $\lim_{|\a| \rightarrow \infty}a_{\a}=0$.

Let us now recall the following  Mahler-type
theorem, a special case of results in \cite{CC, ab0, conrad}:

\begin{theorem}
\label{known}
 Let
$f:\zp \rightarrow \zp$ be a continuous function.  Then there exists a
restricted power series  $F(x_0,x_1,x_2,\ldots )$ in the variables
$x_0,x_1,x_2,\ldots $, with $\zp$-coefficients, such that
$$
f(a)=F(a,\d a, \d^2 a,\ldots )
$$
for all $a$ in $\zp$.
\qed
\end{theorem}

Our Theorems \ref{easy} and \ref{main} imply that the series $F$ in
Theorem \ref{known} can be chosen to depend on finitely many
variables if, and only if, $f$ is analytic.

\subsection{Structure of the paper}
In \S 2 we prove the basic $p$-adic estimates we shall need in our work, and
follow that by proving Theorem \ref{easy} in \S 3. In \S 4 we introduce a
matrix intrinsically associated to the number $p^m$, and analyze its
determinant. They play an important
r\^ole in the proof of Theorem \ref{main}, which we do in \S 5.

\section{$p$-adic estimates}

\begin{lemma}\label{le1}
If $a \in \zp$ and $x=a+p^n u$ in  the disc $a+p^n \zp$, write $\d ^k x$ as a
polynomial function of $u$ of degree $p^k$,
$$\d^k x = \sum_{j=0}^{p^k} c_{a,j}^k u^j \, ,$$
with $c_{a,j}\in \qp$.
Then we have the following  $p$-adic estimates:
\begin{enumerate}
\item $\np{c_{a,0}^k}\leq 1$.
\item $\np{c_{a,1}^k}=\frac{1}{p^{n-k}}.$
\item $\np{c_{a,j}^k}\leq \frac{1}{p^{(n-k+1)j-1}}\, , \quad
 2 \leq j \leq p^k \, .$
\end{enumerate}
\end{lemma}

{\it Proof}. Assertion (1) follows from the equality $c^k_{a,0}=\d^k a$.
We prove assertions (2) and (3) by induction on $k$. The result is
clear for $k=0$. Assuming that it holds for
$k-1$, we prove it for $k$.

By hypothesis we have that
$$\d^{k-1} x =  \sum_{j=0}^{p^{k-1}} c_{a,j}^{k-1} u^j \, ,$$
where for $j\geq 1$ the coefficients $c_{a,j}^{k-1}$ satisfy the estimates
$$
\np{c_{a,j}^{k-1}}\leq \frac{1}{p^{(n-k+2)j-1}}\, ,
$$
with equality for $j=1$.
We use (\ref{do}) to write the $k$-th iterate as
$$
\d^k x = \sum_{j=0}^{p^k} c_{a,j}^k u^j=\frac{ \sum_{j=0}^{p^{k-1}}
c_{a,j}^{k-1} u^j -
( \sum_{j=0}^{p^{k-1}} c_{a,j}^{k-1} u^j )^p }{p} \, .
$$
Consider a fixed index $j\geq 1$. It follows that $c_{a,j}^k$ is equal to
$c^{k-1}_{a,j}/p$ minus a sum of elements of the form a rational
integer times
$$
\frac{c_{a,j_1}^{k-1}\cdots c_{a,j_p}^{k-1}}{p}\, ,
$$
where $j_1+ \cdots +j_{p}=j$. By a commutation, we may assume that
there is an $s$ such that $j_1, \ldots, j_s \geq 1$ and $j_t=0$ for
$t>s$.  Then $s\leq j_1 + \cdots +j_s =j$, and
by assertion (1) of the Lemma and the induction hypothesis applied to each
of the factors, we conclude that
$$
\np{\frac{c_{a,j_1}^{k-1}\cdots c_{a,j_p}^{k-1}}{p}} \leq
\np{\frac{c_{a,j_1}^{k-1}\cdots c_{a,j_s}^{k-1}}{p}}
\leq \frac{1}{p^{(n-k+2)(j_1 + \cdots +j_{s})-s-1}}
\leq \frac{1}{p^{(n-k+1)j-1}}\, ,
$$
and the estimate in (3) follows. The equality in
Assertion (2) follows from the induction hypothesis and the identity
$$
c_{a,1}^k=c_{a,1}^{k-1}\left( \frac{1}{p}-(c_{a,0}^{k-1})^{p-1}\right)\, .
$$
This completes the proof.
\qed

\section{Proof of Theorem \ref{easy}}
Let $f(x)=F(x,\d x, \ldots , \d^m x)$ by an operator of
order $m$ given by the restricted power series $F\in \zp [[t_0, \ldots, t_m]]$.
Thus,
$$
f(x)=\sum_{\a=(\a_0,\ldots, \a_m) }
a_{\a} x^{\a_0}(\d x)^{\a_1} \cdots (\d ^m x)^{\a_m} \, ,
$$
where $a_{\a} \rightarrow 0$ $p$-adically as $|\a|\rightarrow \infty$.

Let $I=\{0,1,\ldots, p^m-1\}$. The family of discs
$\{a+ p^n \zp\}_{a\in I}$ forms a covering of $\zp$.
By Lemma \ref{le1}, if $a \in I$ we have that
$f(a+p^m u)=F_a(u)$, where
$$
F_a(u)=\sum_{\a=(\a_0,\ldots, \a_r) }
a_{\a} (\sum_{j_0=0}^{p^0} c_{a,j_0}^0 u^{j_0})^{\a_0}(\sum_{j_1=0}^{p^1}
c_{a,j_1}^1 u^{j_1})^{\a_1}
 \cdots (\sum_{j_r=0}^{p^m} c_{a,j_m}^m u^{j_m})^{\a_m} \, ,
$$
 with all the $c_{a,j_l}^k$s in $\zp$.
Notice that $F_a(u)$ is a power series in $u$ with $\zp$-coefficients that
go to zero $p$-adically as $|\a|\rightarrow \infty$.
\qed

\section{A matrix associated to $p^m$}
Let us consider the set of all $p$-adic integer roots of the
function  $x \mapsto \d^m x$:
$$
C_m:=\{a \in \zp: \, \d^m a=0\}\, .
$$
Since the $m$-th iterate of $\d$ is given by
a polynomial of degree $p^m$ with $\qp$-coefficients,
$C_m$ has at most $p^m$ elements. In fact, it has exactly $p^m$ elements.
We have:

\begin{lemma}
\label{Cr}
The composition
$$
C_m \subset \zp \rightarrow \zp/p^m\zp
$$
is bijective.
\end{lemma}

{\it Proof}.
We proceed by induction on $m$. For $m=0$ we have that $C_0=\{0\}$. We assume
now that the statement is true for $m-1$, and prove it for $m$.

Given $a \in C_{m-1}$, we consider the polynomial $t^p-t+pa \in \zp[t]$.
By Hensel's lemma, it has $p$ distinct roots that we denote by
$a_1,\ldots ,a_p \in \zp$. Notice that we have $\d a_j=a$ for all $j$s,
and since $\d^{m-1} a=0$, it follows that $\d^m a_j=0$.
We claim that if $a,a' \in C_{m-1}$ and we have that
\begin{equation}
\label{congg}
a_j \equiv a'_{j'}\quad \text{mod $p^m$}\, ,
\end{equation}
for some $j, j'$, then $a=a'$ and $j=j'$. Indeed,
if (\ref{congg}) holds then $a \equiv a'$ mod $p^{m-1}$, and by
the induction hypothesis, $a=a'$ and hence $j=j'$ as well.
By this claim and the induction hypothesis, $C_m$ contains a set of $p^m$
elements that injects into $\zp/p^m\zp$. As $C_m$ has at most
$p^m$ elements, this forces the map $C_m\rightarrow \zp/p^m\zp$ to be
bijective.
\qed

By Lemma \ref{Cr}, we can write
\begin{equation}
\label{CCr}
C_m=\{a_0,a_1, \ldots , a_{p^m-1}\}\, ,
\end{equation}
where $a_{\alpha} \equiv \alpha$ mod $p^m$ for all $\alpha$s in the set
\begin{equation}
\label{AA}
I=\{0,\ldots ,p^m -1\}\, .
\end{equation}
On the other hand, let us fix an ordering in the set
\begin{equation}
I'=\{\beta=(\beta_0,\ldots ,\beta_{m-1})\in \z^m: \, 0 \leq \beta_0,\ldots ,
\beta_{m-1} \leq p-1\}\, ,
\label{BB}
\end{equation}
and consider the $(p^m -1)\times (p^m -1)$-matrix
\begin{equation}
W=(w_{\alpha \beta})_{\a\in I, \b \in I'}\, ,
\label{WWVV}
\end{equation}
whose entries are given by
$$
w_{\alpha \beta}:=(a_{\alpha})^{\beta_0} (\d a_{\alpha})^{\beta_1}
\ldots (\d^{m-1} a_{\alpha})^{\beta_{m-1}}\in \zp\, ,
$$
where we have used the convention that $a^0=1$ for all $a \in \zp$.
 Up to a permutation of its columns,
the matrix $W$ is intrinsically associated to the number $p^m$.

\begin{lemma}
\label{VW}
The determinant of the matrix $W$
is invertible in $\zp$.
\end{lemma}

{\it Proof}. We use the reduction mod $p$ mapping
$$
\begin{array}{ccl}
\zp & \rightarrow & {\mathbb F}_p:=\zp/p\zp \\
a & \mapsto & \overline{a}
\end{array}\, .
$$
By Lemma 3.20 in \cite{ab2}, the function
$$
\begin{array}{ccl}
\zp & \rightarrow & {\mathbb F}_{p}^m \\
a & \mapsto & (\overline{a},\overline{\d a},\ldots , \overline{\d^{m-1} a})
\end{array}
$$
induces a bijection
$\zp/p^m\zp \simeq {\mathbb F}_p^m$. For any
$\gamma=(\gamma_0,\ldots ,\gamma_{m-1}) \in {\mathbb F}_p^m$
and any $\beta \in I'$, we set
$$
v_{\gamma \beta}=\gamma_0^{\beta_0} \gamma_1^{\beta_1} \cdots
\gamma_{m-1}^{\beta_{m-1}}\, .
$$
Notice that $\d^i a_{\alpha} \equiv \d^i \alpha$ mod $p$ if $i \leq m-1$.
Therefore, by Lemma \ref{Cr}, in order to prove  Lemma \ref{VW} we just
need to show that $\det{(v_{\gamma \beta})}\neq 0\in {\mathbb F}_p$.

Assume the latter is false. This means that there exist constants
$\lambda_{\beta_0 \ldots \beta_{m-1}} \in {\mathbb F}_p$
for $(\beta_0,\ldots ,\beta_{m-1}) \in I'$, not all zero, such that
$$
\sum_{\beta_0=0}^{p-1} \ldots
\sum_{\beta_{m-1}=0}^{p-1}\lambda_{\beta_0 \ldots \beta_{m-1}}
\gamma_0^{\beta_0} \gamma_1^{\beta_1} \ldots \gamma_{m-1}^{\beta_{m-1}}=0
$$
for all $\gamma\in {\mathbb F}_p^m$. By induction on $m$, this easily implies
that all the $\lambda$s vanish, a contradiction. This proves our Lemma.
\qed

\section{Proof of Theorem \ref{main}}
We now carry out the proof of Theorem \ref{main} by proving the
following result that is more precise.
In what follows, $I$ and $I'$ are the sets of indices
(\ref{AA}) and (\ref{BB}), respectively.

\begin{theorem}
\label{mainn}
Let $f:\zp \rightarrow\zp$ be an analytic function of level $m$.
Then there exists a unique restricted power series
$F \in \zp[[x_0,x_1,\ldots ,x_m]]$ with the following properties:
\begin{enumerate}
\item $F(x_0,x_1,\ldots ,x_m)=\sum_{n\geq 0}
\sum_{\b \in I'} a_{\b,n}x_0^{\b_0}x_1^{\b_1}\cdots
x_{m-1}^{\b_{m-1}} x_m^{n}$.
\item $f(a)=F(a,\d a,\ldots ,\d^m a)$, $a \in \zp$.
\end{enumerate}
\end{theorem}

{\it Proof}. We start by proving the existence of $F$.
Notice that if $g:\zp\rightarrow \zp$ is any arithmetic differential
operator of order $m$ and $a \in \zp$, then $h(x):=g(x+a)$ is also an
arithmetic differential operator of order $m$; cf. Remark \ref{trei}, (3).
By this, and without losing generality, we may assume
that the function $f$ in the statement of the Theorem is zero on all discs
of radius $1/p^m$ except for $p^m\zp$. Without losing generality also, we
may additionally assume that there exists an
 $l\geq 0$ such that $f(p^mu)=u^l$ for all $u \in \zp$.

We recall the set $C_m$ in (\ref{CCr}).
The family of discs $\{a+p^m \zp\}_{a\in C_m}$ forms
a covering of $\zp$. Notice that $a_0=0$.

By Lemma \ref{le1}, for  $a \in C_m$ and $0\leq k\leq m$, we have
that $\d^k(a+p^mu)=\sum_{j=0}^{p^m} c_{a,j}^k u^j$, with
$c_{a,0}^k=\d^k a$, and
\begin{equation}
\label{zuzu}
\np{c_{a,0}^k}\leq 1\, , \quad \np{c_{a,1}^k}=\frac{1}{p^{m-k}}\, ,
\quad \np{c_{a,j}^k} \leq \frac{1}{p^{(m-k+1)j-1}}\, , \;
2\leq j \leq p^k\, .
\end{equation}
We may view $\d^k(a+p^mu)$ as an element in the ring of polynomials $\zp[u]$.
Since $\d^m a=0$, it follows that $c^m_{a,j}=0$.

We proceed to inductively determine the coefficients $a_{\b,n}$ in the
series $F$ appearing in the statement of the Theorem, so that
\begin{equation}
\np{a_{\b,n}}\leq \min\left\{ 1, \frac{1}{p^{n-l}}\right\}
 \, , \quad n\geq 0\, , \; \b \in I'\, ,
\label{est}
\end{equation}
and so that,
 if $F_{a}(u)=F(a+p^m u)$, for $a \in C_m$, then we have
\begin{equation}
\left\{ \begin{array}{lclll}
F_a(u)& = & u^l & \text{if} & a=0\, , \\ F_a(u) & = & 0 & \text{if} & a \neq 0.
\end{array}\right. \label{eqp}
\end{equation}
Here we view (\ref{eqp}) as equalities of functions of $u \in \zp$.
However, since each $F_a(u)$ is defined by a restricted power series in
$\zp[[u]]$, it is enough to check (\ref{eqp}) as equalities
in the ring of formal power series $\zp[[u]]$.

We consider the polynomials $F_{a}^k(u) \in \zp[u]$ defined by
$$
F_{a}^k(u):=
\sum_{n=0}^k \sum_{\b \in I'}
a_{\beta,n}\left(\sum_{j=0}^{p^0} c^0_{a,j}u^j\right)^{\beta_0}\ldots
\left(\sum_{j=0}^{p^{m-1}} c^{m-1}_{a,j}u^j\right)^{\beta_{m-1}}
\left(\sum_{j=1}^{p^m} c^m_{a,j}u^j \right)^n
$$
so that $F_a^k(u)$ converge $u$-adically to $F_a$ in $\zp[[u]]$.
We find the $a_{\b,n}$s inductively so they satisfy the estimate
(\ref{est}), and such that the following congruences hold in
the ring $\zp[u]$:
\begin{equation}
\left\{ \begin{array}{lclllll}
F^k_a(u)& \equiv & u^l & \text{mod} & u^{k+1} & \text{if} & a=0\, , \\
F^k_a(u) & \equiv & 0 & \text{mod} & u^{k+1} & \text{if} & a \neq 0.
\end{array}\right. \label{scallo}
\end{equation}

In what follows we denote by $\d_{ij}$ the Kronecker symbol.

For the starting point of the induction, we choose the coefficients
$a_{\b,0}$, $\b \in I'$, such that (\ref{scallo}) and (\ref{est}) hold.
This can be achieved by solving the
system of equations
$$
\sum_{\beta \in I'} w_{\alpha \beta}
a_{\beta,0}=\d_{l0}\d_{\alpha 0}\, , \quad \a \in I\, ,
$$
where $W=(w_{\alpha \beta})$ is the  matrix (\ref{WWVV}).
By Lemma \ref{VW}, this system
can be readily solved for the $a_{\beta,0}$s, with
the solution being a vector of $p$-adic integers.

For the $k$-th step of the induction, let us notice that for $a=a_{\a}$,
the coefficient of $u^k$ in $F_a^k(u)$  is
given by
\begin{equation}
\label{coeef}
(c_{a,1}^m)^k
\sum_{\beta \in I'} w_{\alpha\beta} a_{\beta,k}+
\sum_{n=0}^{k-1}\sum_{\beta \in I'}a_{\beta,n} b_{\beta,n,k}\, ,
\end{equation}
where $b_{\beta,n,k}$ is the coefficient of $u^k$ in
$$
\left(\sum_{j=0}^{p^0} c^0_{a,j}u^j\right)^{\beta_0}\ldots
\left(\sum_{j=0}^{p^{m-1}} c^{m-1}_{a,j}u^j\right)^{\beta_{m-1}}
\left(\sum_{j=1}^{p^m} c^m_{a,j}u^j \right)^n \, .
$$
Thus,
$b_{\beta,n,k}$ is a $\z$-linear combination of products of the
form
$$
\left(\prod_{i=1}^{\beta_0}c^0_{a,j_{0i}}\right)\ldots
\left(\prod_{i=1}^{\beta_{m-1}} c^{m-1}_{a,j_{m-1,i}}\right)
\left(\prod_{i=1}^n
c^m_{a,j_{mi}}\right)\, ,
$$
with
$$\sum_{r=0}^{m-1}\sum_{i=1}^{\beta_r}j_{ri}+
\sum_{i=1}^n j_{mi}=k\, .
$$
We may assume that
there are integers $s_r$ such that $j_{ri}\geq 1$ for $i\leq s_r$
and $j_{ri}=0$ for $i>s_r$. So we have
\begin{equation}
\label{xix}
s_r\leq \sum_{i=1}^{s_r} j_{ri}\, ,\quad
\sum_{r=0}^{m-1} \sum_{i=1}^{s_r}j_{ri}+
\sum_{i=1}^n j_{mi}=k\, .
\end{equation}
By (\ref{zuzu}) and the induction hypothesis,
\begin{equation}
\np{a_{\beta,n} b_{\beta,n,k}} \leq
\min\left\{1,\frac{1}{p^{n-l+\sigma}}\right\}\, .
\end{equation}
where
$$
\sigma=\sum_{r=0}^{m-1}[(m-r+1)(\sum_{i=1}^{s_r} j_{ri}) -s_r]+\sum_{i=1}^n
j_{mi}-n\, .
$$
Now, by (\ref{xix}) we have that
$$
\begin{array}{rcl}
\sigma & \geq &  \sum_{r=0}^{m-1}
[2 (\sum_{i=1}^{s_r} j_{ri})-s_r ]+\sum_{i=1}^n j_{mi}-n \vspace{2mm} \\
 & \geq & \sum_{r=0}^{m-1} \sum_{i=1}^{s_r} j_{ri}+
\sum_{i=1}^n j_{mi}-n \vspace{2mm} \\
 & = & k-n\, .
\end{array}
$$
Hence
\begin{equation}
\label{lappo}
\np{a_{\beta,n} b_{\beta,n,k}} \leq \min\left\{1,\frac{1}{p^{k-l}}\right\}\, .
\end{equation}
Now, by the induction hypothesis also, $F_a^k(u)$ satisfies (\ref{scallo})
if we have
\begin{equation}\label{eqk}
\sum_{\beta \in I'} w_{\alpha \beta} a_{\beta,k} =(c^m_{a_{\alpha},1})^{-k}
\left(\d_{kl}\d_{\alpha 0}-\sum_{n=0}^{k-1}
\sum_{\beta \in I'}a_{\beta,n} b_{\beta,n,k}\right)\, , \; \alpha \in I.
\end{equation}
By (\ref{zuzu}) and (\ref{lappo}),
the $p$-adic norm of the right hand side of (\ref{eqk})
is bounded above by
$\min\{ 1, 1/p^{k-l}\}$.
Again, by Lemma \ref{WWVV}, we can solve the system
(\ref{eqk}) for the $a_{\beta,k}$s, with the solution satisfying the
estimates (\ref{est}).
This completes the induction, and hence the existence part of the Theorem.

In order to prove the uniqueness, we need to show that if a restricted power
series $F$ satisfies conditions (1) and (2) in the Theorem for $f=0$,
then $a_{\beta,n}=0$ for all $\beta\in I'$, $n \geq 0$. This follows
by an induction on $n$, in view of the equalities
$$
\sum_{\beta \in I'} w_{\alpha \beta} a_{\beta,k} =-(c^m_{a_{\alpha},1})^{-k}
\left(\sum_{n=0}^{k-1}
\sum_{\beta \in I'}a_{\beta,n} b_{\beta,n,k}\right)\, , \; \alpha \in I.
$$
This finishes the proof.
\qed

\end{document}